\documentclass{article}

\usepackage{graphicx}
\newtheorem{thm}{Theorem}[section]
\newtheorem{cor}[thm]{Corollary}
\newtheorem{lem}[thm]{Lemma}
\newtheorem{de}[thm]{Definition}
\newcounter{bean}
\newcounter{milk}

\begin{document}

\title{The Mystery of the Shape Parameter III}         % Enter your title between curly braces
\author{Lin-Tian Luh\\Department of Financial and Computational Mathematics\\ Providence University\\
Shalu Area, Taichung City, Taiwan\\ Email:ltluh@pu.edu.tw}        % Enter your name between curly braces
\date{\today}

          % Enter your date or \today between curly braces
\maketitle

{\bf Abstract.} This is a continuation of our earlier study of the shape parameter c contained in the famous multiquadrics $(-1)^{\lceil \frac{\beta}{2} \rceil}(c^{2}+\|x\|^{2})^{\frac{\beta}{2}},\ \beta>0$, and the inverse multiquadrics $(c^{2}+\|x\|^{2})^{\beta},\ \beta <0$. In the previous two papers we presented criteria for the optimal choice of c, based on the exponential-type error bound. In this paper a new set of criteria is developed, based on the improved exponential-type error bound. This results in much sharper error estimates when c is chosen appropriately, with the same size of fill distance. What is important is that the optimal value of c can be successfully predicted without any search when fill distance is of reasonable size, making it practically useful. The drawback is that the distribution of the data points is not purely scattered. However it seems to be harmless.

\section{Introduction}       % Enter section title between curly braces
We begin with some basic ingredients of our theoretical ground.

Let $T_{n}$ denote the $n$-simplex in $R^{n}$ whose definition can be found in \cite{Fl}. A $1$-simplex is a line segment, a $2$-simplex is a triangle, and a $3$-simplex is a tetrahedron with four vertices.

Let $v_{i},\ 1\leq i\leq n+1$ be the vertices of $T_{n}$. Then any point $x\in T_{n}$ can be written as a convex combination of the vertices:
$$x=\sum_{i=1}^{n+1}c_{i}v_{i},\ \sum_{i=1}^{n+1}c_{i}=1,\ c_{i}\geq 0.$$
The numbers $c_{1},\cdots, c_{n+1}$ are called the barycentric coordinates of $x$.

For any $n$-simplex $T_{n}$, the evenly spaced points of degree $l$ are those points whose barycentric coordinates are of the form
$$(\frac{k_{1}}{l},\frac{k_{2}}{l}, \cdots, \frac{k_{n+1}}{l}),\ k_{i}\ nonnegative\ integers\ with \ \sum_{i=1}^{n+1}k_{i}=l.$$ If we let $P_{l}^{n}$ denote the space of polynomials of degree not exceeding $l$ in n variables, it is easily seen that the number of evenly spaced points of degree $l$ is exactly $N=\dim P_{l}^{n}=\left( \begin{array}{c}
                                                         n+l \\
                                                          n
                                                       \end{array} \right) $. Also, such points form a determining set for $P_{l}^{n}$, by \cite{Bo}.

In this paper the interpolation will happen in an $n$-simplex and the set $X$ of centers(interpolation points) will be evenly spaced points in the $n$-simplex.

The radial function we use is
\begin{eqnarray}
  h(x):=\Gamma(-\frac{\beta}{2})(c^{2}+\|x\|^{2})^{\frac{\beta}{2}},\ \beta\in R\backslash 2N_{\geq 0},\ c>0,
\end{eqnarray}
where $\|x\|$ is the Euclidean norm of $x\in R^{n}$, $\Gamma$ is the classical gamma function, and $\beta,\ c$ are constants. Note that this definition is slightly different from the one mentioned in the abstract. We adopt (1) because it will greatly simplify its Fourier transform and our future work. The function $h(x)$ in (1) is conditionally positive definite(c.p.d.) of order $m=\max \{0,\lceil \frac{\beta}{2}\rceil \}$ where $\lceil \frac{\beta}{2}\rceil $ means the smallest integer greater than or equal to $\frac{\beta}{2}$. Further details can be found in \cite{MN1, We}.

Given data points $(x_{j},y_{j}),\ j=1,\cdots,N$, where $X=\{ x_{1},\cdots,x_{N}\}$ is a subset of $R^{n}$ and $y_{j}$ are real or complex numbers, our interpolant will be of the form
\begin{eqnarray}
  s(x)=p(x)+\sum_{j=1}^{N}c_{j}h(x-x_{j}),
\end{eqnarray}
where $p(x)$ is a polynomial in $P_{m-1}^{n}$ to be determined and $c_{j}$ are coefficients to be chosen.

As is well known in the theory of radial basis functions, if $X$ is a determining set for $P_{m-1}^{n}$, there exists a unique polynomial $p(x)$ and unique constants $c_{1},\cdots,c_{N}$ satisfying the linear system
\begin{eqnarray}
  p(x_{i})+\sum_{j=1}^{N}c_{j}h(x_{i}-x_{j})&=&y_{i},\ \ i=1,\cdots,N \nonumber\\
                                                             \\
                 \sum_{j=1}^{N}c_{j}q(x_{j})&=&0 \nonumber
\end{eqnarray}
where $q$ ranges over all basis elements of $P_{m-1}^{n}$. All these can be found in \cite{MN1, We}.

\subsection{Fundamental theory}    % Enter subsection title between curly braces
Each function $h$ of the form (1) induces a function space ${\cal C}_{h,m}$, called native space, whose definition and characterization can be found in \cite{MN1,MN2,Lu1,Lu2,Lu4-1,We}. Here $m=\max \{0,\lceil \frac{\beta}{2} \rceil \}$. Also, there is a seminorm $\|f\|_{h}$ for each $f\in {\cal C}_{h,m}$. In our theory every interpolated function belongs to the native space.

Before entering the main theorem, let us introduce two constants.
\begin{de}
Let $n$ and $\beta$ be as in (1). The numbers $\rho$ and $\Delta_{0}$ are defined as follows.
\begin{list}
  {(\alph{bean})}{\usecounter{bean} \setlength{\rightmargin}{\leftmargin}}
  \item Suppose $\beta <n-3$. Let $s=\lceil \frac{n-\beta -3}{2}\rceil $. Then 
    \begin{list}{(\roman{milk})}{\usecounter{milk} \setlength{\rightmargin}{\leftmargin}}
      \item if $\beta <0,\ \rho=\frac{3+s}{3}\ and\  \Delta_{0}=\frac{(2+s)(1+s)\cdots 3}{
 \rho^{2}};$
      \item if $\beta >0,\ \rho=1+\frac{s}{2\lceil \frac{\beta}{2}\rceil +3} \ and \ \Delta_{0}=\frac{(2m+2+s)(2m+1+s)\cdots (2m+3)}{\rho^{2m+2}}$ \\
where $ m=\lceil \frac{\beta}{2}\rceil$.          
    \end{list}
  \item Suppose $n-3\leq \beta <n-1$. Then $\rho=1$ and $\Delta_{0}=1$.
  \item Suppose $\beta \geq n-1$. Let $s=-\lceil \frac{n-\beta -3}{2}\rceil $. Then
 $$\rho =1\ and \ \Delta_{0}=\frac{1}{(2m+2)(2m+1)\cdots (2m-s+3)} \ where \ m=\left \lceil \frac{\beta}{2}\right \rceil.$$  
\end{list}
\end{de}
Our criteria for the optimal choice of c is based on the following theorem which we take directly from \cite{Lu4} but with a slight modification to make it easier to understand.
\begin{thm}
  Let $h$ be as in (1). For any positive number $b_{0}$, let $C=\max \left\{ \frac{2}{3b_{0}},\frac{8\rho}{c}\right\}$ and $\delta_{0}=\frac{1}{3C}$. For any n-simplex $Q$ of diameter $r$ satisfying $\frac{1}{3C}\leq r\leq \frac{2}{3C}$(note that $\frac{2}{3C}\leq b_{0}$), if $f\in {\cal C}_{h,m}$,
\begin{eqnarray}
  |f(x)-s(x)|\leq 2^{\frac{n+\beta-7}{4}}\pi^{\frac{n-1}{4}}\sqrt{n\alpha_{n}}c^{\frac{\beta}{2}}\sqrt{\Delta_{0}}\sqrt{3C}\sqrt{\delta}(\lambda')^{\frac{1}{\delta}}\|f\|_{h}
\end{eqnarray}
holds for all $x\in Q$ and $0<\delta\leq\delta_{0}$, where $s(x)$ is defined as in (2) with $x_{1},\cdots ,x_{N}$ the evenly spaced points of degree $l$ in $Q$ satisfying $\frac{1}{3C\delta}\leq l\leq \frac{2}{3C\delta}$. The constant $\alpha_{n}$ denotes the volume of the unit ball in $R^{n}$, and $0<\lambda'<1$ is given by 
$$\lambda'=\left(\frac{2}{3}\right)^{\frac{1}{3C}}$$
which only in some cases mildly depends on the dimension n.

\end{thm}
{\bf Remark}: This seemingly complicated theorem is in fact not difficult to understand. Note that the right-hand side of (4) approaches zero as $\delta$ tends to zero. Hence $\delta$ is in spirit like the well-known fill distance, although not exactly the same. Also, the upper bound in (4) is greatly influenced by the shape parameter $c$. The only thing which is not transparent is the relation between $c$ and $\|f\|_{h}$. Consequently, in order to make it useful in the choice of $c$, we still have to do some work.

We begin with the following definition.
\begin{de}
  For any $\sigma>0$, the class of band-limited functions $f$ in $L^{2}(R^{n})$ is defined by
$$B_{\sigma}=\{f\in L^{2}(R^{n}):\ \hat{f}(\xi)=0\ if\ |\xi|>\sigma\},$$ where $\hat{f}$ denotes the Fourier transform of $f$.
\end{de}
Now we cite Theorem1.6 of \cite{Lu5} as a lemma.
\begin{lem}
  Let $h$ be as in (1) with $\beta>0$. Any function $f$ in $B_{\sigma}$ belongs to ${\cal C}_{h,m}$ and 
$$\|f\|_{h}\leq \sqrt{m!S(m,n)}2^{-n-\frac{1+\beta}{4}}\pi^{-n-\frac{1}{4}}\sigma^{\frac{1+\beta+n}{4}}e^{\frac{c\sigma}{2}}c^{\frac{1-\beta-n}{4}}\|f\|_{L^{2}(R^{n})}$$
where $c,\beta$ are as in (1) and $S(m,n)$ is a constant determined by $m$ and $n$.
\end{lem}

\begin{cor}
  Let $h$ be as in (1) with $\beta>0$. If $\sigma>0$ and $f\in B_{\sigma}$, the inequality (4) can be transformed into
\begin{eqnarray}
  |f(x)-s(x)|\leq 2^{-2-\frac{3}{4}n}\pi^{\frac{-2-3n}{4}}\sqrt{n\alpha_{n}}\sqrt{\Delta_{0}}\sqrt{3C}\sqrt{m!S(m,n)}\sigma^{\frac{1+\beta+n}{4}}c^{\frac{1+\beta-n}{4}}e^{\frac{c\sigma}{2}}\sqrt{\delta}(\lambda')^{\frac{1}{\delta}}\|f\|_{L^{2}(R^{n})}
\end{eqnarray}
\end{cor}
In order to handle the case $\beta<0$, we need the following lemma which is just Theorem1.7 of \cite{Lu5}.
\begin{lem}
  Let $h$ be as in (1) with $\beta<0$ such that $n+\beta\geq 1$ or $n+\beta=-1$. Any function $f$ in $B_{\sigma}$ belongs to ${\cal C}_{h,m}$ and satisfies
$$\|f\|_{h}\leq 2^{-n-\frac{1+\beta}{4}}\pi^{-n-\frac{1}{4}}\sigma^{\frac{1+\beta+n}{4}}e^{\frac{c\sigma}{2}}c^{\frac{1-n-\beta}{4}}\|f\|_{L^{2}(R^{n})}.$$
\end{lem}

\begin{cor}
  Let $h$ be as in (1) with $\beta<0$ such that $n+\beta\geq 1$ or $n+\beta=-1$. If $\sigma>0$ and $f\in B_{\sigma}$, the inequality (4) can be transformed into
\begin{eqnarray}
|f(x)-s(x)|\leq 2^{-2-\frac{3}{4}n}\pi^{\frac{-3n-2}{4}}\sqrt{n\alpha_{n}}\sqrt{\Delta_{0}}\sqrt{3C}\sigma^{\frac{1+\beta+n}{4}}c^{\frac{1+\beta-n}{4}}e^{\frac{c\sigma}{2}}\sqrt{\delta}(\lambda')^{\frac{1}{\delta}}\|f\|_{L^{2}(R^{n})}
\end{eqnarray}
\end{cor}
Note that Corollary1.7 does not cover the frequently seen case $\beta=-1,n=1$. For this we need the following lemma which is just Lemma2.1 of \cite{Lu5}.

\begin{lem}
  Let $h$ be as in (1) with $\beta=-1,n=1$. For any $\sigma>0$, if $f\in B_{\sigma}$, then $f\in {\cal C}_{h,m}$ and 
$$\|f\|_{h}\leq (2\pi)^{-n}2^{-\frac{1}{4}}\left\{ \frac{1}{{\cal K}_{0}(1)}\int_{|\xi|\leq \frac{1}{c}}|\hat{f}(\xi)|^{2}d\xi +\frac{1}{a_{0}}\int_{\frac{1}{c}<|\xi|\leq \sigma}|\hat{f}(\xi)|^{2}\sqrt{c|\xi|}e^{c|\xi|}d\xi \right\} ^{1/2}$$
if $\frac{1}{c}<\sigma$, where $a_{0}=\frac{1}{2\sqrt{3}}$, and
$$\|f\|_{h}\leq (2\pi)^{-n}2^{-\frac{1}{4}}\left\{ \frac{1}{{\cal K}_{0}(1)}\int_{|\xi|\leq \frac{1}{c}}|\hat{f}(\xi)|^{2}d\xi \right\}^{1/2}$$
if $\frac{1}{c}\geq \sigma$.
\end{lem}

\begin{cor}
  Let $h$ be as in (1) with $\beta=-1,n=1$. If $\sigma>0$ and $f\in B_{\sigma}$, the inequality (4) can be transformed into 
\begin{eqnarray}
  |f(x)-s(x)|\leq 2^{-2+\frac{-3n+\beta}{4}}\pi^{\frac{-3n-1}{4}}\sqrt{n\alpha_{n}}\sqrt{\Delta_{0}}\sqrt{3C}c^{\frac{\beta}{2}}\sqrt{\delta}(\lambda')^{\frac{1}{\delta}}\{A+B\}^{1/2}
\end{eqnarray}
where $A=\frac{1}{{\cal K}_{0}(1)}\int_{|\xi|\leq \frac{1}{c}}|\hat{f}(\xi)|^{2}d\xi$ for all c, $B=2\sqrt{3}\int_{\frac{1}{c}<|\xi|\leq \sigma}|\hat{f}(\xi)|^{2}\sqrt{c|\xi|}e^{c|\xi|}d\xi$ if $\frac{1}{c}<\sigma$, and $B=0$ if $\frac{1}{c}\geq \sigma$.
\end{cor}

\section{Criteria of Choosing c}
The results of Section1 provide us with useful theoretical ground for choosing c. Note that in the right-hand side of (5),(6) and (7), there is always a main function determined by c. Let us call it the MN function, denoted by $MN(c)$, as in \cite{Lu6}. Its graph is called the MN curve. Then finding the optimal value of c is equivalent to finding the minimum of $MN(c)$. The range of $c$ should be clarified first. In order to satisfy the condition $\delta \leq \delta_{0}$ as required by Theorem1.2, for any given $b_{0}$, we require that $\delta< \frac{b_{0}}{2}$ and $c\in [24\rho \delta, \infty)$. Then we have three cases as follows. \\
\\
 {\bf Case1}. \fbox{$\beta>0$nd $n\geq 1$} Let $f\in B_{\sigma}$ and $h$ be defined as in (1) with $\beta>0$ and $n\geq1$. For any given $b_{0}>0$ and $\delta<\frac{b_{0}}{2}$, under the conditions of Theorem1.2, the optimal choice of $c$ in the interval $[24\rho \delta, \infty)$ is the number minimizing 
$$MN(c):= \left\{ \begin{array}{ll}
                    \sqrt{8\rho}\cdot c^{\frac{\beta-1-n}{4}}\cdot e^{\frac{c\sigma}{2}}\cdot\left(\frac{2}{3}\right) ^{\frac{c}{24\delta\rho}}   & \mbox{if $24\rho \delta \leq c<12\rho b_{0},$} \\
                    \sqrt{\frac{2}{3b_{0}}}\cdot c^{\frac{1+\beta-n}{4}}\cdot e^{\frac{c\sigma}{2}}\cdot \left( \frac{2}{3}\right) ^{\frac{b_{0}}{2\delta}} 
                            & \mbox{if $c\geq 12\rho b_{0}$.}

                  \end{array}
          \right.  $$        
{\bf Reason}: This follows directly from (5).\\
\\
{\bf Examples}:\\
\\

\begin{figure}[t]
\centering
\includegraphics[scale=1.0]{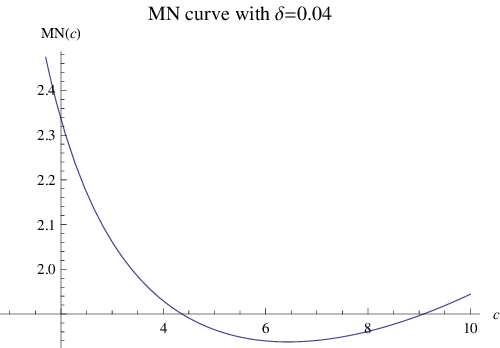}
\caption{Here $n=2,\ \beta=1, b_{0}=10$ and $\sigma=1$.}

\includegraphics[scale=1.0]{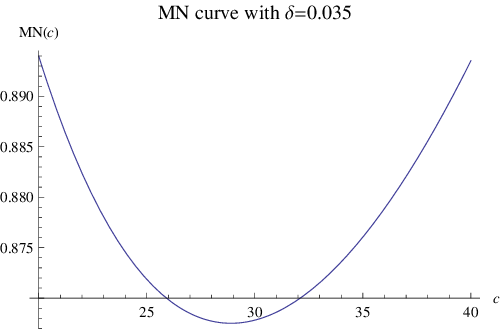}
\caption{Here $n=2,\ \beta=1, b_{0}=10$ and $\sigma=1$.}

\includegraphics[scale=1.0]{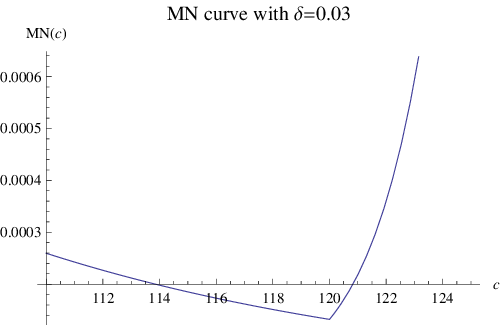}
\caption{Here $n=2,\ \beta=1, b_{0}=10$ and $\sigma=1$.}

\end{figure}

\clearpage
\begin{figure}[t]
\centering
\includegraphics[scale=1.0]{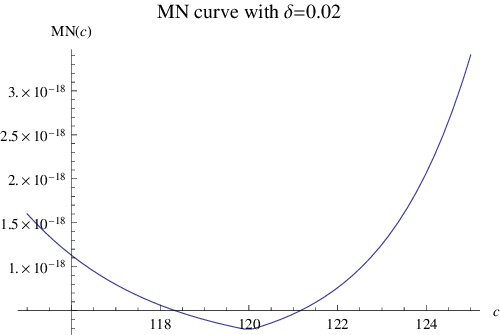}
\caption{Here $n=2,\ \beta=1, b_{0}=10$ and $\sigma=1$.}

\includegraphics[scale=1.0]{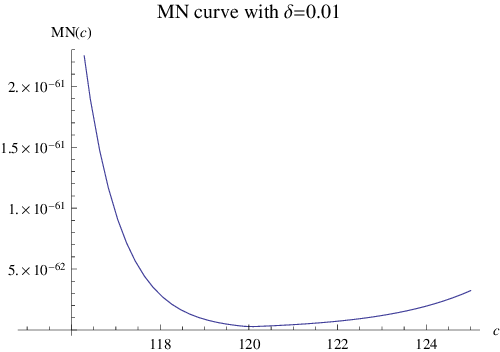}
\caption{Here $n=2,\ \beta=1, b_{0}=10$ and $\sigma=1$.}

\end{figure}
For $\beta<0$, we separate it into two cases.\\
\\
{\bf Case2}. \fbox{$\beta<0$ and $n+\beta\geq1$ or $n+\beta=-1$} Let $f\in B_{\sigma}$ and $h$ be defined as in (1) with $\beta<0$ and $n+\beta\geq 1$, or $n+\beta=-1$. For any given $b_{0}>0$ and $\delta<\frac{b_{0}}{2}$, under the conditions of Theorem1.2, the optimal choice of $c$ in the interval $[24\rho \delta, \infty)$ is the number minimizing
$$MN(c):= \left\{ \begin{array}{ll}
                    \sqrt{8\rho}\cdot c^{\frac{\beta-1-n}{4}}\cdot e^{\frac{c\sigma}{2}}\cdot\left(\frac{2}{3}\right) ^{\frac{c}{24\delta\rho}}   & \mbox{if $24\rho \delta \leq c<12\rho b_{0},$} \\
                    \sqrt{\frac{2}{3b_{0}}}\cdot c^{\frac{1+\beta-n}{4}}\cdot e^{\frac{c\sigma}{2}}\cdot \left( \frac{2}{3}\right) ^{\frac{b_{0}}{2\delta}} 
                            & \mbox{if $c\geq 12\rho b_{0}$.}

                  \end{array}
          \right.  $$
{\bf Reason}: This is an immediate result of Corollary1.7. \\
\\
{\bf Examples}:\\
\\

\begin{figure}[t]
\centering
\includegraphics[scale=1.0]{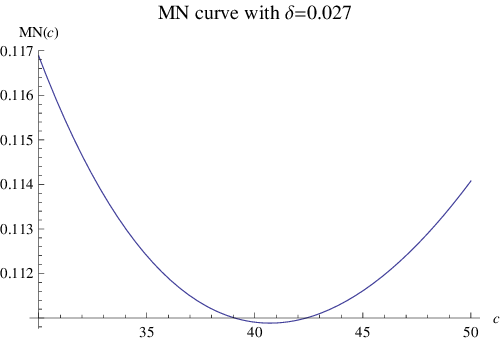}
\caption{Here $n=3,\ \beta=-1, b_{0}=20$ and $\sigma=1$.}

\includegraphics[scale=1.0]{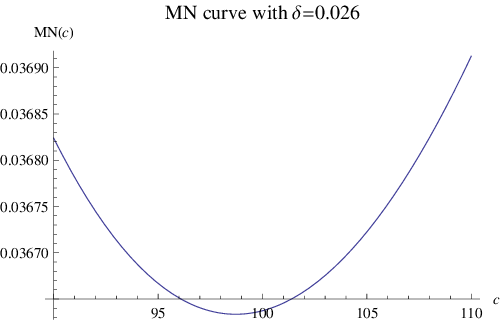}
\caption{Here $n=3,\ \beta=-1, b_{0}=20$ and $\sigma=1$.}

\includegraphics[scale=1.0]{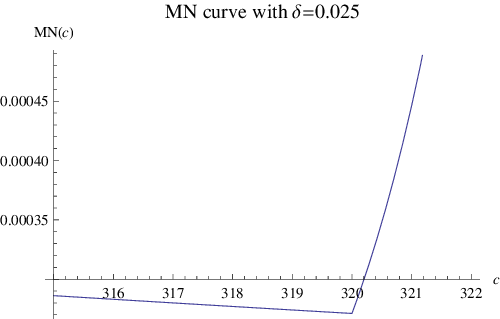}
\caption{Here $n=3,\ \beta=-1, b_{0}=20$ and $\sigma=1$.}

\end{figure}

\clearpage

\begin{figure}[t]
\centering
\includegraphics[scale=1.0]{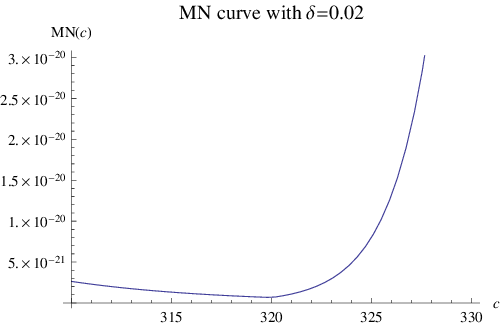}
\caption{Here $n=3,\ \beta=-1, b_{0}=20$ and $\sigma=1$.}

\includegraphics[scale=1.0]{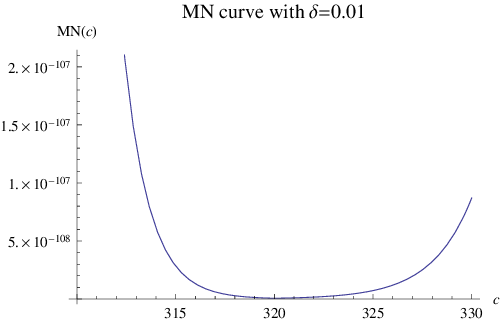}
\caption{Here $n=3,\ \beta=-1, b_{0}=20$ and $\sigma=1$.}

\end{figure}
Now we begin the case $\beta=-1$ and $n=1$.\\
\\
{\bf Case3}. \fbox{$\beta=-1$ and $n=1$} Let $ f\in B_{\sigma}$ and $h$ be defined as in (1) with $\beta=-1$ and $n=1$. For any given $b_{0}>0$ and $\delta<\frac{b_{0}}{2}$, under the conditions of Theorem1.2, the optimal choice of $c$ in the interval $[24\rho \delta, \infty)$ is the number minimizing

$$MN(c):=\left\{ \begin{array}{ll}
                   \sqrt{8\rho}\cdot c^{\frac{\beta-1}{2}}\cdot \left( \frac{2}{3}\right) ^{\frac{c}{24\delta\rho}}M(c)   & \mbox{if $24\rho \delta\leq c<12\rho b_{0},$} \\
                   \sqrt{\frac{2}{3b_{0}}}\cdot c^{\frac{\beta}{2}}\cdot \left( \frac{2}{3}\right) ^{\frac{b_{0}}{2\delta}}M(c)& \mbox{if $c\geq 12\rho b_{0},$} 
                \end{array}
        \right. $$ 
where $M(c)$ is defined by
$$M(c):=\left\{  \begin{array}{ll}
                  \frac{1}{\sqrt{{\cal K}_{0}(1)}}   & \mbox{if $c\leq \frac{1}{\sigma},$} \\
                  \left\{ \frac{1}{{\cal K}_{0}(1)}+2\sqrt{3}\sqrt{c\sigma}e^{c\sigma}\right\}^{1/2} & \mbox{if $c>\frac{1}{\sigma},$}
                \end{array}
        \right. $$
${\cal K}_{0}$ being the modified Bessel function.\\
\\
{\bf Reason}: Note that in (7), $\{A+B\}^{1/2}$ can be further treated as follows. 
For $0<c\leq \frac{1}{\sigma}$, we have $\sigma\leq \frac{1}{c}$ and
$$A=\frac{1}{{\cal K}_{0}(1)}\int_{|\xi|\leq\frac{1}{c}}|\hat{f}(\xi)|^{2}d\xi=\frac{1}{{\cal K}_{0}(1)}\int_{|\xi|\leq \sigma}|\hat{f}(\xi)|^{2}d\xi=\frac{1}{{\cal K}_{0}(1)}\|f\|^{2}_{L^{2}(R^{n})}$$ because $f\in B_{\sigma}$. Therefore
$$\{A+B\}^{1/2}=\sqrt{A}=\frac{1}{\sqrt{{\cal K}_{0}(1)}}\|f\|_{L^{2}(R^{n})}$$ if $0<c\leq \frac{1}{\sigma}$. Now, if $\frac{1}{\sigma}<c$,
\begin{eqnarray*}
  \{A+B\}^{1/2} & = & \left\{ \frac{1}{{\cal K}_{0}(1)}\int_{|\xi|\leq \frac{1}{c}}|\hat{f}(\xi)|^{2}d\xi+2\sqrt{3}\int_{\frac{1}{c}<|\xi|\leq \sigma}|\hat{f}(\xi)|^{2}\sqrt{c|\xi|}e^{c|\xi|}d\xi \right\}^{1/2}\\
                                             & \leq & \left\{ \frac{1}{{\cal K}_{0}(1)}\int_{|\xi|<\sigma}|\hat{f}(\xi)|^{2}d\xi+2\sqrt{3}\int_{|\xi|\leq \sigma}|\hat{f}(\xi)|^{2}\sqrt{c|\xi|}e^{c|\xi|}d\xi\right\}^{1/2}\\
                                             & \leq & \left\{ \frac{1}{{\cal K}_{0}(1)}\|f\|^{2}_{L^{2}(R^{n})}+2\sqrt{3}\sqrt{c\sigma}e^{c\sigma}\|f\|^{2}_{L^{2}(R^{n})}\right\}^{1/2}\\
                                             & = & \left\{ \frac{1}{{\cal K}_{0}(1)}+2\sqrt{3}\sqrt{c\sigma}e^{c\sigma}\right\}^{1/2}\|f\|_{L^{2}(R^{n})}.
\end{eqnarray*}
Our conclusion thus follows.\\
\\
{\bf Examples}:\\
\\
\begin{figure}[h]
\centering
\includegraphics[scale=1.0]{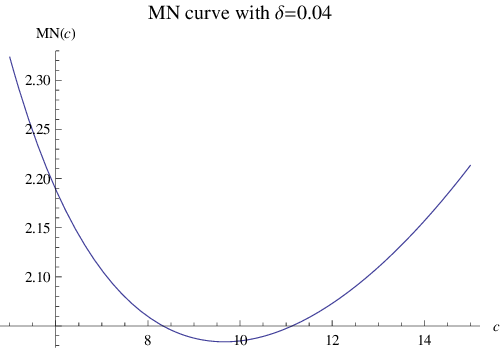}
\caption{Here $n=1,\ \beta=-1,\ b_{0}=5$ and $\sigma=1$.} 
\end{figure}

\begin{figure}[t]
\centering
\includegraphics[scale=1.0]{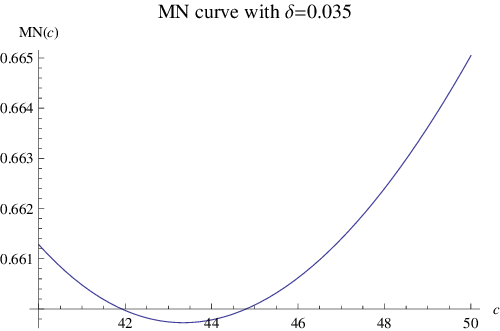}
\caption{Here $n=1,\ \beta=-1,\ b_{0}=5$ and $\sigma=1$.}

\includegraphics[scale=1.0]{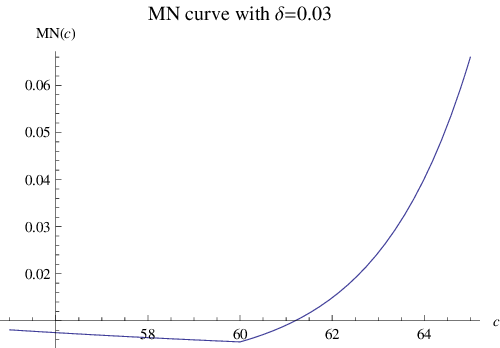}
\caption{Here $n=1,\ \beta=-1,\ b_{0}=5$ and $\sigma=1$.}

\includegraphics[scale=1.0]{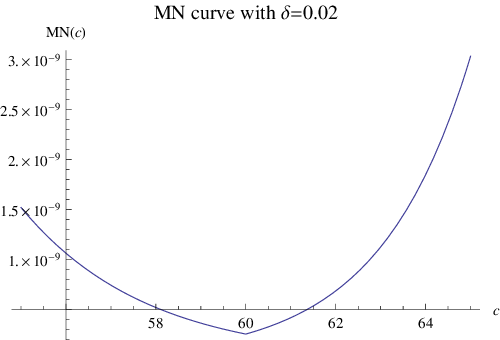}
\caption{Here $n=1,\ \beta=-1,\ b_{0}=5$ and $\sigma=1$.}

\end{figure}

\clearpage

\begin{figure}[t]
\centering
\includegraphics[scale=1.0]{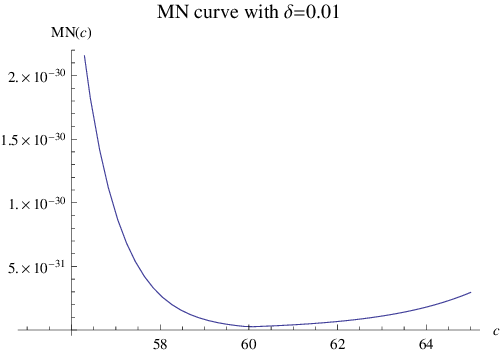}
\caption{Here $n=1,\ \beta=-1,\ b_{0}=5$ and $\sigma=1$.}

\end{figure}
Note that the optimal c increases rapidly as $\delta$ becomes small.

\section{Experiment}
In this section we test Case1 of the preceding section and let $n=2,\ \beta=1,\ b_{0}=10,\ \sigma=1$. In order to make it more useful and understandable, we replace the function $h(x)$ defined in Section1 by the more commonly used function $h(x):=-(c^{2}+\|x\|^{2})^{\frac{1}{2}},\ x\in R^{2}$. The interpolation occurs in a regular triangle with side length $r$. As required by Theorem1.2, $1/(3C)\leq r\leq 2/(3C)$. We choose to let $r=2/(3C)$ so that the centers(interpolation points) will not be too close to each other. By the very definition of $C$, $C$ depends on $c$. Therefore, as $c$ changes, the diameter $r$ of the triangle also changes. Let the original vertices be $v_{1}=(-1,0),\ v_{2}=(1,0)$, and $v_{3}=(0,\sqrt{3})$. Then the triangle we adopt has vertices $w_{1}=\frac{r}{2}v_{1}, w_{2}=\frac{r}{2}v_{2}$, and $w_{3}=\frac{r}{2}v_{3}$. As a result, the side length will be $c/12$ if $c<120$, and $10$ if $c\geq120$. The centers are evenly spaced ponits of degree $l$ in the triangle. Theorem1.2 requires that $1/(3C\delta)\leq l\leq 2/(3C\delta)$. The smaller $l$ is, the less data points will be used, as can be seen in the beginning of Section1. Hence we choose $l=1/(3C\delta)$. Once the centers are arranged, we let the test points be the evenly spaced points in the same triangle with degree $l'=l+1$. We adopt the root-mean-square error to evaluate the distance between the approximated and approximating functions at the test points. Let $f$ and $s$ denote the approximated and approximating functions, respectively. If the test points are $z_{1},\cdots , z_{N}$, then $$RMS=\sqrt{\sum_{i=1}^{N} (f(z_{i})-s(z_{i}))^{2}/N}$$ is its root-mean-square error. We denote the number of data points and the number of test points by $N_{d}$ and $N_{t}$, respectively. Thus, if the centers are evenly spaced points of degree $l$, then $N_{d}=(l+2)(l+1)/2$, and $N_{t}=(l+3)(l+2)/2$, respectively.

There is a crucial logical problem in using our approach. Note that in our core theorem Theorem1.2, the function $h$, and hence the shape parameter $c$, appears first. Theoretically, one should fix $c$ and then choose the other parameters, including $\delta$, which is in spirit like the well-known fill distance. However, since the optimal $c$ cannot be known in advance, we choose the other parameters first. Then choose the optimal $c$ according to the MN curve. Once $c$ is chosen, we begin to design the simplex(triangle in the $R^{2}$ case) and the interpolation points(centers) in the simplex according to Theorem1.2. This is just a trick for avoiding logical troubles.

In this experiment, the approximated function is $$f(x,y):=\frac{\sin{(\frac{\sigma x}{\sqrt{2}})}}{\frac{\sigma x}{\sqrt{2}}}\cdot \frac{\sin{(\frac{\sigma y}{\sqrt{2}})}}{\frac{\sigma y}{\sqrt{2}}}$$ where we let $\frac{\sin{z}}{z}=1$ if $z=0$. The map $f$ can be easily checked to belong to $B_{\sigma}$ defined in Definition1.3.

We emphasize that our approach for choosing $c$ optimally is reliable only when the parameter $\delta$ is small enough. This phenomenon can also be seen in the experiment of \cite{Lu8}. Note that in Figures 1-5, the value $c$ minimizing the MN curve moves rapidly to 120 and remains there when $\delta$ is small. It strongly suggests that one should choose $c=120$ as the optimal value when $\delta$ is small enough. We test $\delta=0.225,\ 0.2,\ 0.175,\ 0.15,\ 0.125,\ 0.1,\ 0.085,$ and $0.075$. 

In the following tables, we use $COND$, $N_{d}$ and $N_{t}$ to denote the condition number of the interpolation matrix, the number of data points used, and the number of test points, respectively. The condition number is the traditional one, i.e., the infinity-norm condition number. In virtue of the arbitrarily precise computer software Mathematica, the problem of ill-conditioning is resolved by adopting enough effective digits to the right of the decimal point for each calculation when the condition number is very large, at the cost of spending considerable computer time. 

\begin{table}[h]
\caption{$\delta=0.225$}
\centering
\tiny
\begin{tabular}{c llllll}\\[2ex]
\hline\hline \\ [1ex]
\large $c$ & \large $80$ & \large $100$ & \large $110$ & \large $115$ & \large $117$ & \large $118$  \\ [1ex]
\hline\\[1ex]

\large $RMS$ & \normalsize $ 2.0\cdot 10^{-13} $ & \normalsize $8.8\cdot 10^{-16} $ & \normalsize $ 6.1\cdot 10^{-17} $ & \normalsize $ 3.9\cdot 10^{-17} $ & \normalsize $6.0\cdot 10^{-17} $  & \normalsize $7.4\cdot 10^{-17}$   \\ [1ex]
\hline\\[1ex]

\large $COND$ & \normalsize $7.4\cdot 10^{53}$ & \normalsize $1.0\cdot 10^{68}$ & \normalsize $1.3 \cdot 10^{75}$ & \normalsize $4.1 \cdot 10^{78}$ & \normalsize $4.1 \cdot 10^{78}$ & \normalsize $4.1 \cdot 10^{78}$ \\ [1ex]
\hline\\[1ex]

\large $N_{d}$ & \normalsize $136$  & \normalsize $210$  & \normalsize $253$  & \normalsize $276$ & \normalsize $276$  & \normalsize $276$\\[1ex]
\hline\\[1ex]

\large $N_{t}$ & \normalsize $153$  & \normalsize $231$  & \normalsize $276$  & \normalsize $300$ & \normalsize $300$  & \normalsize $300$  \\[1ex]

\hline\hline \\[1ex]

\large $c$ & \large $119$ & \large $120$ & \large $130$ & \large $140$ & \large $160$  \\ [1ex]
\hline\\[1ex]

\large $RMS$ & \normalsize $ 3.5\cdot 10^{-18} $ & \normalsize $4.3\cdot 10^{-18} $ & \normalsize $ 4.8\cdot 10^{-18} $ & \normalsize $ 5.3\cdot 10^{-18} $ & \normalsize $6.0\cdot 10^{-18} $     \\ [1ex]
\hline \\[1ex]

\large $COND$ & \normalsize $1.4\cdot 10^{82}$ & \normalsize $1.4\cdot 10^{82}$ & \normalsize $5.5 \cdot 10^{83}$ & \normalsize $1.7 \cdot 10^{85}$ & \normalsize $7.8 \cdot 10^{87}$  \\ [1ex]
\hline \\[1ex]

\large $N_{d}$ & \normalsize $300$   & \normalsize $300$  & \normalsize $300$ & \normalsize $300$  & \normalsize $300$\\[1ex]
\hline \\[1ex]

\large $N_{t}$ & \normalsize $325$  & \normalsize $325$  & \normalsize $325$ & \normalsize $325$  & \normalsize $325$ \\[1ex]
\hline\\[1ex]

\end{tabular}
\label{001}
\end{table}

\clearpage

\begin{table}[t]
\caption{$\delta=0.2$}
\centering
\tiny
\begin{tabular}{c llllll}\\[2ex]
\hline\hline \\ [1ex]
\large $c$ & \large $80$ & \large $100$ & \large $110$ & \large $115$ & \large $116$ & \large $117$  \\ [1ex]
\hline\\[1ex]

\large $RMS$ & \normalsize $ 1.5\cdot 10^{-15} $ & \normalsize $7.8\cdot 10^{-18} $ & \normalsize $ 5.8\cdot 10^{-19} $ & \normalsize $ 3.6\cdot 10^{-19} $ & \normalsize $1.8\cdot 10^{-20} $  & \normalsize $2.2\cdot 10^{-20}$   \\ [1ex]
\hline\\[1ex]

\large $COND$ & \normalsize $8.5\cdot 10^{60}$ & \normalsize $1.3\cdot 10^{75}$ & \normalsize $1.4 \cdot 10^{82}$ & \normalsize $5.0 \cdot 10^{85}$ & \normalsize $1.6 \cdot 10^{89}$ & \normalsize $1.6 \cdot 10^{89}$ \\ [1ex]
\hline\\[1ex]

\large $N_{d}$ & \normalsize $171$  & \normalsize $253$  & \normalsize $300$  & \normalsize $325$ & \normalsize $351$  & \normalsize $351$\\[1ex]
\hline\\[1ex]

\large $N_{t}$ & \normalsize $190$  & \normalsize $276$  & \normalsize $325$  & \normalsize $351$ & \normalsize $378$  & \normalsize $378$ \\[1ex]

\hline\hline \\[1ex]

\large $c$ & \large $118$ & \large $119$ & \large $120$ & \large $130$ & \large $140$  & \large $160$\\ [1ex]
\hline\\[1ex]

\large $RMS$ & \normalsize $ 2.8\cdot 10^{-20} $ & \normalsize $3.4\cdot 10^{-20} $ & \normalsize $ 4.2\cdot 10^{-20} $ & \normalsize $ 5.0\cdot 10^{-20} $ & \normalsize $5.6\cdot 10^{-20} $  & \normalsize $6.7\cdot 10^{-20}$   \\ [1ex]
\hline \\[1ex]

\large $COND$ & \normalsize $1.6\cdot 10^{89}$ & \normalsize $1.6\cdot 10^{89}$ & \normalsize $1.6 \cdot 10^{89}$ & \normalsize $8.8 \cdot 10^{90}$ & \normalsize  $3.6 \cdot 10^{92}$ & \normalsize $2.8\cdot 10^{95}$ \\ [1ex]
\hline \\[1ex]

\large $N_{d}$ & \normalsize $351$  & \normalsize $351$  & \normalsize $351$  & \normalsize $351$ & \normalsize $351$  & \normalsize $351$ \\[1ex]
\hline \\[2ex]

\large $N_{t}$ & \normalsize $378$  & \normalsize $378$  & \normalsize $378$  & \normalsize $378$ & \normalsize $378$  & \normalsize $378$ \\[1ex]
\hline\\[1ex]

\end{tabular}
\label{002}
\end{table}

\begin{table}[h]
\caption{$\delta=0.175$}
\centering
\tiny
\begin{tabular}{c llllll}\\[2ex]
\hline\hline \\ [2ex]
\large $c$ & \large $80$ & \large $100$ & \large $110$ & \large $116$ & \large $117$ & \large $118$  \\ [1ex]
\hline\\[1ex]

\large $RMS$ & \normalsize $ 6.1\cdot 10^{-19} $ & \normalsize $6.8\cdot 10^{-21} $ & \normalsize $ 2.9\cdot 10^{-23} $ & \normalsize $ 2.5\cdot 10^{-23} $ & \normalsize $3.3\cdot 10^{-23} $  & \normalsize $1.5\cdot 10^{-24}$   \\ [1ex]
\hline\\[1ex]

\large $COND$ & \normalsize $3.5\cdot 10^{71}$ & \normalsize $5.0\cdot 10^{85}$ & \normalsize $1.9 \cdot 10^{96}$ & \normalsize $6.2 \cdot 10^{99}$ & \normalsize $6.2 \cdot 10^{99}$ & \normalsize $2.1 \cdot 10^{103}$ \\ [1ex]
\hline\\[2ex]

\large $N_{d}$ & \normalsize $231$ & \normalsize $325$  & \normalsize $406$ & \normalsize $435$  & \normalsize $435$ & \normalsize $465$ \\[1ex]
\hline\\[1ex]

\large $N_{t}$ & \normalsize $253$  & \normalsize $351$ & \normalsize $435$ & \normalsize $465$ & \normalsize $465$ & \normalsize $496$ \\[1ex]

\hline\hline \\[1ex]

\large $c$ & \large $119$ & \large $120$ & \large $130$ & \large $140$ & \large $160$  \\ [1ex]
\hline\\[1ex]

\large $RMS$ & \normalsize $ 1.9\cdot 10^{-24} $ & \normalsize $2.5\cdot 10^{-24} $ & \normalsize $ 3.3\cdot 10^{-24} $ & \normalsize $ 4.1\cdot 10^{-24} $ & \normalsize $5.5\cdot 10^{-24} $    \\ [1ex]
\hline \\[1ex]

\large $COND$ & \normalsize $2.1\cdot 10^{103}$ & \normalsize $2.1\cdot 10^{103}$ & \normalsize $2.2 \cdot 10^{105}$ & \normalsize $1.6 \cdot 10^{107}$ & \normalsize  $3.7 \cdot 10^{110}$  \\ [1ex]
\hline \\[1ex]

\large $N_{d}$ & \normalsize $465$  & \normalsize $465$  & \normalsize $465$ & \normalsize $465$  & \normalsize $465$  \\[1ex]
\hline \\[1ex]

\large $N_{t}$ & \normalsize $496$  & \normalsize $496$  & \normalsize $496$ & \normalsize $496$  & \normalsize $496$ \\[1ex]
\hline\\[1ex]

\end{tabular}
\label{003}
\end{table}

\clearpage

\begin{table}[t]
\caption{$\delta=0.15$}
\centering
\tiny
\begin{tabular}{c llllll}\\[2ex]
\hline\hline \\ [1ex]
\large $c$ & \large $80$ & \large $100$ & \large $110$ & \large $118$ & \large $119$ & \large $120$  \\ [1ex]
\hline\\[1ex]

\large $RMS$ & \normalsize $ 1.3\cdot 10^{-22} $ & \normalsize $1.3\cdot 10^{-25} $ & \normalsize $ 6.7\cdot 10^{-28} $ & \normalsize $ 3.8\cdot 10^{-29} $ & \normalsize $6.6\cdot 10^{-30} $  & \normalsize $9.6\cdot 10^{-30}$   \\ [1ex]
\hline\\[1ex]

\large $COND$ & \normalsize $1.4\cdot 10^{82}$ & \normalsize $6.2\cdot 10^{99}$ & \normalsize $2.4 \cdot 10^{110}$ & \normalsize $2.7 \cdot 10^{117}$ & \normalsize $8.8 \cdot 10^{120}$ & \normalsize $8.8 \cdot 10^{120}$ \\ [1ex]
\hline\\[2ex]

\large $N_{d}$ & \normalsize $300$  & \normalsize $435$  & \normalsize $528$ & \normalsize $595$  & \normalsize $630$ & \normalsize $630$\\[1ex]
\hline\\[1ex]

\large $N_{t}$ & \normalsize $325$  & \normalsize $465$ & \normalsize $561$ & \normalsize $630$ & \normalsize $666$ & \normalsize $666$ \\[1ex]

\hline\hline \\[1ex]

\large $c$ & \large $121$ & \large $122$ & \large $130$ & \large $140$ & \large $160$  \\ [1ex]
\hline\\[1ex]

\large $RMS$ & \normalsize $ 1.0\cdot 10^{-29} $ & \normalsize $1.1\cdot 10^{-29} $ & \normalsize $ 1.8\cdot 10^{-29} $ & \normalsize $ 2.9\cdot 10^{-29} $ & \normalsize $5.3\cdot 10^{-29} $    \\ [1ex]
\hline \\[1ex]

\large $COND$ & \normalsize $1.5\cdot 10^{121}$ & \normalsize $2.7\cdot 10^{121}$ & \normalsize $2.0 \cdot 10^{123}$ & \normalsize $3.1 \cdot 10^{125}$ & \normalsize  $2.7 \cdot 10^{129}$  \\ [1ex]
\hline \\[1ex]

\large $N_{d}$ & \normalsize $630$   & \normalsize $630$  & \normalsize $630$ & \normalsize $630$  & \normalsize $630$ \\[1ex]
\hline \\[1ex]

\large $N_{t}$ & \normalsize $666$  & \normalsize $666$  & \normalsize $666$ & \normalsize $666$  & \normalsize $666$  \\[1ex]
\hline\\[1ex]

\end{tabular}
\label{004}
\end{table}

\begin{table}[h]
\caption{$\delta=0.125$}
\centering
\tiny
\begin{tabular}{c llllll}\\[2ex]
\hline\hline \\ [1ex]
\large $c$ & \large $80$ & \large $100$ & \large $110$ & \large $114$ & \large $116$ & \large $117$  \\ [1ex]
\hline\\[1ex]

\large $RMS$ & \normalsize $ 4.2\cdot 10^{-28} $ & \normalsize $7.1\cdot 10^{-34} $ & \normalsize $ 8.3\cdot 10^{-36} $ & \normalsize $ 2.8\cdot 10^{-36} $ & \normalsize $2.5\cdot 10^{-37} $  & \normalsize $4.3\cdot 10^{-37}$   \\ [1ex]
\hline\\[1ex]

\large $COND$ & \normalsize $1.9\cdot 10^{96}$ & \normalsize $8.8\cdot 10^{120}$ & \normalsize $3.2 \cdot 10^{131}$ & \normalsize $1.1 \cdot 10^{135}$ & \normalsize $3.7 \cdot 10^{138}$ & \normalsize $3.7 \cdot 10^{138}$ \\ [1ex]
\hline\\[1ex]

\large $N_{d}$ & \normalsize $406$  & \normalsize $630$ & \normalsize $741$ & \normalsize $780$ & \normalsize $820$ & \normalsize $820$\\[1ex]
\hline\\[1ex]

\large $N_{t}$ & \normalsize $435$  & \normalsize $666$ & \normalsize $780$ & \normalsize $820$ & \normalsize $861$ & \normalsize $861$ \\[1ex]

\hline\hline \\[1ex]

\large $c$ & \large $118$ & \large $119$ & \large $120$ & \large $130$ & \large $140$  & \large $160$\\ [1ex]
\hline\\[1ex]

\large $RMS$ & \normalsize $ 2.4\cdot 10^{-38} $ & \normalsize $4.8\cdot 10^{-38} $ & \normalsize $ 9.1\cdot 10^{-38} $ & \normalsize $ 5.1\cdot 10^{-37} $ & \normalsize $1.4\cdot 10^{-36} $  & \normalsize $4.4\cdot 10^{-36}$   \\ [1ex]
\hline \\[1ex]

\large $COND$ & \normalsize $1.2\cdot 10^{142}$ & \normalsize $1.2\cdot 10^{142}$ & \normalsize $1.2 \cdot 10^{142}$ & \normalsize $7.1 \cdot 10^{144}$ & \normalsize  $2.7 \cdot 10^{147}$  & \normalsize $1.2\cdot 10^{152}$\\ [1ex]
\hline \\[1ex]

\large $N_{d}$ & \normalsize $861$  & \normalsize $861$  & \normalsize $861$ & \normalsize $861$  & \normalsize $861$ & \normalsize $861$  \\[1ex]
\hline \\[1ex]

\large $N_{t}$ & \normalsize $903$  & \normalsize $903$  & \normalsize $903$ & \normalsize $903$  & \normalsize $903$ & \normalsize $903$   \\[1ex]
\hline\\[1ex]

\end{tabular}
\label{005}
\end{table}

\clearpage

\begin{table}[t]
\caption{$\delta=0.1$}
\centering
\tiny
\begin{tabular}{c llllll}\\[2ex]
\hline\hline \\ [1ex]
\large $c$ & \large $80$ & \large $100$ & \large $110$ & \large $116$ & \large $117$ & \large $118$  \\ [1ex]
\hline\\[1ex]

\large $RMS$ & \normalsize $ 3.9\cdot 10^{-28} $ & \normalsize $7.2\cdot 10^{-45} $ & \normalsize $ 2.8\cdot 10^{-48} $ & \normalsize $ 3.0\cdot 10^{-51} $ & \normalsize $5.2\cdot 10^{-51} $  & \normalsize $1.4\cdot 10^{-52}$   \\ [1ex]
\hline\\[1ex]

\large $COND$ & \normalsize $8.8\cdot 10^{120}$ & \normalsize $1.3\cdot 10^{149}$ & \normalsize $1.5 \cdot 10^{163}$ & \normalsize $5.5 \cdot 10^{173}$ & \normalsize $5.5 \cdot 10^{173}$ & \normalsize $1.8 \cdot 10^{177}$ \\ [1ex]
\hline\\[1ex]

\large $N_{d}$ & \normalsize $630$ & \normalsize $946$ & \normalsize $1128$ & \normalsize $1275$ & \normalsize $1275$ & \normalsize $1326$\\[1ex]
\hline\\[1ex]

\large $N_{t}$ & \normalsize $666$  & \normalsize $990$ & \normalsize $1176$ & \normalsize $1326$ & \normalsize $1326$ & \normalsize $1378$\\[1ex]

\hline\hline \\[1ex]

\large $c$ & \large $119$ & \large $120$ & \large $130$ & \large $140$ & \large $160$  \\ [1ex]
\hline\\[1ex]

\large $RMS$ & \normalsize $ 4.4\cdot 10^{-52} $ & \normalsize $9.9\cdot 10^{-52} $ & \normalsize $ 7.7\cdot 10^{-52} $ & \normalsize $ 7.7\cdot 10^{-51} $ & \normalsize $2.6\cdot 10^{-49} $     \\ [1ex]
\hline \\[1ex]

\large $COND$ & \normalsize $1.8\cdot 10^{177}$ & \normalsize $1.8\cdot 10^{177}$ & \normalsize $5.4 \cdot 10^{180}$ & \normalsize $8.9 \cdot 10^{183}$ & \normalsize  $5.6 \cdot 10^{189}$  \\ [1ex]
\hline \\[1ex]

\large $N_{d}$ & \normalsize $1326$  & \normalsize $1326$ & \normalsize $1326$  & \normalsize $1326$ & \normalsize $1326$ \\[1ex]
\hline \\[1ex]

\large $N_{t}$ & \normalsize $1378$  & \normalsize $1378$ & \normalsize $1378$  & \normalsize $1378$ & \normalsize $1378$ \\[1ex]
\hline\\[1ex]

\end{tabular}
\label{006}
\end{table}

Note that when the parameter $\delta$ decreases, the condition numbers get large. For $\delta=0.1$ and $c=120$, we adopted 200 effective digits to the right of the decimal point for each calculation and successfully overcame the problem of ill-conditioning. The other cases were handled in a similar way. As we emphasized, our approach of choosing $c$ optimally is reliable only when $\delta$ is small enough. We of course want to decrease $\delta$ further until the optimal $c$ coincides with the theoretical value completely. However, limited by the speed of the computer, we have to reduce the scale of our experiment for $\delta<0.1$. In the following two tables, we only test five values of $c$ for each $\delta$.

\begin{table}[h]
\caption{$\delta=0.085$}
\centering
\tiny
\begin{tabular}{c lllll}\\[2ex]
\hline\hline \\ [1ex]
\large $c$ & \large $110$ & \large $118$ & \large $119$ & \large $120$ & \large $130$  \\ [1ex]
\hline\\[1ex]

\large $RMS$ & \normalsize $ 1.6\cdot 10^{-59} $ & \normalsize $1.7\cdot 10^{-63} $ & \normalsize $ 2.6\cdot 10^{-65} $ & \normalsize $ 1.4\cdot 10^{-64} $ & \normalsize $1.6\cdot 10^{-64} $     \\ [1ex]
\hline\\[1ex]

\large $COND$ & \normalsize $2.2\cdot 10^{191}$ & \normalsize $2.4\cdot 10^{205}$ & \normalsize $8.0 \cdot 10^{208}$ & \normalsize $8.0 \cdot 10^{208}$ & \normalsize $1.0 \cdot 10^{213}$  \\ [1ex]
\hline\\[1ex]

\large $N_{d}$ & \normalsize $1540$ & \normalsize $1770$ & \normalsize $1830$ & \normalsize $1830$  & \normalsize $1830$\\[1ex]
\hline\\[1ex]

\large $N_{t}$ & \normalsize $1596$  & \normalsize $1830$ & \normalsize $1891$ & \normalsize $1891$  & \normalsize $1891$        \\[1ex]
\hline\\[1ex]

\end{tabular}
\label{007}
\end{table}

\clearpage

\begin{table}[t]
\caption{$\delta=0.075$}
\centering
\tiny
\begin{tabular}{c lllll}\\[2ex]
\hline\hline \\ [1ex]
\large $c$ & \large $110$ & \large $115$ & \large $118$ & \large $119$ & \large $120$   \\ [1ex]
\hline\\[1ex]

\large $RMS$ & \normalsize $2.5\cdot 10^{-70}$ & \normalsize $ 5.6\cdot 10^{-72} $ & \normalsize $5.4\cdot 10^{-74} $ & \normalsize $ 3.0\cdot 10^{-75} $ & \normalsize $ 4.2\cdot 10^{-75} $    \\ [1ex]
\hline\\[1ex]

\large $COND$ & \normalsize $2.8\cdot 10^{219}$ & \normalsize $3.0\cdot 10^{226}$ & \normalsize $3.2\cdot 10^{233}$ & \normalsize $1.0 \cdot 10^{237}$ & \normalsize $1.0 \cdot 10^{237}$   \\ [1ex]
\hline\\[1ex]

\large $N_{d}$ & \normalsize $2016$ & \normalsize $2145$ & \normalsize $2278$ & \normalsize $2346$ & \normalsize $2346$\\[1ex]
\hline\\[1ex]

\large $N_{t}$ & \normalsize $2080$   & \normalsize $2211$  & \normalsize $2346$ & \normalsize $2415$  & \normalsize $2415$\\[1ex]
\hline\\[1ex]

\end{tabular}
\label{008}
\end{table}

It can be seen that in these tables, the optimal $c$ tends to be moving to the theoretical value 120 as $\delta$ decreases. Among them, the case $\delta=0.075$ is most important because $\delta$ is smaller. To our regret, there is still a very small gap between the experimentally optimal value and the theoretically predicted one. We have reason to believe that if $\delta$ is further decreased, they will coincide completely, as can be seen in the experiment of \cite{Lu8}. In this paper we cannot do so because it takes too much computer time. Even for $\delta=0.075$, it requires two and half hours to complete only one command. In order to test one $c$, at least five hours must be spent. If $\delta$ is further decreased, maybe 30 hours will be needed to test only one value of $c$.

As a whole, these results are quite satisfactory and our approach of choosing the shape parameter can be trusted.

\section{Summary}
Both \cite{Lu5} and this paper deal with the interpolation of band-limited functions. In \cite{Lu5} the range of c is $[12\rho\sqrt{n}e^{2n\gamma_{n}}\gamma_{n}(m+1)\delta, \infty)$ where $\delta$ is the well-known fill distance and $\gamma_{n}$'s are integers which grow very fast as the dimension n increases. In order to make the left endpoint of the closed-open interval small enough, often $\delta$ must be very small. The consequence is that a huge number of data points will be involved, making the criteria of choosing c only theoretically valuable, especially for $n\geq 2$. Now we have greatly improved this restriction and replaced the left endpoint by $24\rho \delta$ which is much smaller. In fact, we can further enlarge the range of c and allow $c\in (0,\infty)$. However, the interpolation domain(the simplex) required by our main theorem will become very small and is not worth doing. Experiments show that our criteria apply even for c less than $24\rho\delta$. We just cannot prove it. Consequently the restriction $c\geq 24\rho\delta$ does not seem to be a big problem. 

What is important is that our criteria of choosing $c$ are based on the error bound presented in Theorem1.2. After all, error bound and error are not exactly the same. What we can control is error bound, not error. Maybe the choice of $c$ is also influenced by the closeness of the shapes of the approximated and the approximating functions. If the two surfaces match each other well, the root-mean-square error will be small even when the error bound is not small. In any case, empirical results show that our criteria are very reliable. Even if there is a gap between the experimental and theoretical values, the gap is very small.

As for the function space, although it is required that the approximated function should belong to the $B_{\sigma}$ space, our approach in fact applies to any function in the Sobolev space. As shown in \cite{NWW,We2}, any function in the Sobolev space can be interpolated by a $B_{\sigma}$ function with a good error bound. Then any $B_{\sigma}$ function can be interpolated by MQ(multiquadrics) or IMQ(inverse multiquadrics), also with a good error bound. Thus the function in the Sobolev space can be interpolated by MQ and IMQ, with the same set of data points. The $B_{\sigma}$ function plays only an intermediate role and need not be found explicitly. One needs only to know its existence. The error estimate can be handled by triangle inequality. In other words, when dealing with functions in the Sobolev space, we already know how to choose the shape parameter contained in MQ and IMQ. This is particularly meaningful in solving partial differential equations with MQ and IMQ because a lot of important PDE's have solutions in the Sobolev space.

% Set the ending of a LaTeX document
\end{document}